\documentclass[12pt]{article}
\newcommand{\di}{\displaystyle}
\usepackage{amsmath}
\usepackage{amssymb}
\usepackage{amsthm}
\newtheorem{theorem}{Theorem}

\newtheorem{definition}{Definition}
\newtheorem{corollary}{Corollary}
\newtheorem{remark}{Remark}
\begin{document}
\title{ Some Operators Associated to Rarita-Schwinger Type Operators}
\author{Junxia Li and John Ryan\\
\emph{\small Department of Mathematics, University of Arkansas, Fayetteville, AR 72701, USA}}
\date{}

\maketitle
\begin{abstract}
In this paper we study some operators associated to the Rarita-Schwinger operators. They arise from the difference between the Dirac operator and the Rarita-Schwinger operators. These operators are called remaining operators. They are based on the Dirac operator and projection operators $I-P_k.$ The fundamental solutions of these operators are harmonic polynomials, homogeneous of degree $k$. First we study the remaining operators and their representation theory in Euclidean space. Second, we can extend the remaining operators in Euclidean space to the sphere under the Cayley transformation.
\end{abstract}

\\
\\
\par \emph{This paper is dedicated to Michael Shapiro on the occasion of his 65th birthday.}

\section{Introduction}
 Rarita-Schwinger operators in Clifford analysis arise in representation theory for the covering groups of $SO(n)$ and $O(n)$. They are generalizations of the Dirac operator. We denote a Rarita-Schwinger operator by $R_k$, where $k=0, 1, \cdots, m, \cdots.$ When $k=0$ it is the Dirac operator. The Rarita-Schwinger operators $R_k$ in Euclidean space have been studied in \cite{BSSV, BSSV1,DLRV, Va1, Va2}.  Rarita-Schwinger operators on the sphere denoted by $R_k^S$ have also been studied in \cite{LRV}.

 In this paper we study the remaining operators, $Q_k$, which are related to the Rarita-Schwinger operators. In fact, The remaining operators are the difference between the Dirac operator and the Rarita-Schwinger operators.

Let $\mathcal{H}$$_{k}$ be the space of harmonic polynomials homogeneous of degree $k$ and $\mathcal{M}_{k}$, $\mathcal{M}_{k-1}$ be the spaces of $Cl_n-$ valued monogenic polynomials, homogeneous of degree $k$ and $k-1$ respectively. Instead of considering $P_k:  \mathcal{H}_k\rightarrow \mathcal{M}_{k}$ in \cite{DLRV}, we look at the projection map, $I-P_k:  \mathcal{H}_k\rightarrow u\mathcal{M}_{k-1},$ and the Dirac operator, to define the $Q_k$ operators and construct their fundamental solutions in $\mathbb{R}^n$. We introduce basic results of these operators. This includes Stokes' Theorem, Borel-Pompeiu Theorem, Cauchy's Integral Formula and a Cauchy Transform. In section 5, by considering the Cayley transformation and its inverse, we can extend the results for the remaining operators in $\mathbb{R}^n$ to the sphere, $\mathbb{S}^n.$ We construct the fundamental solutions to the remaining operators by applying the Cayley transformation to the fundamental solutions to the remaining operators in $\mathbb{R}^n.$ We also obtain the intertwining operators for $Q_k$ operators and $Q_k^S$ operators. In turn, we establish the conformal invariance of the remaining equations under the Cayley transformation and its inverse. We conclude by giving some basic integral formulas with detailed proofs, and pointing out that results obtained for Rarita-Schwinger operators in \cite{LRV} for real projective space readily carry over to the context presented here.

\section{Preliminaries}
\par A Clifford algebra, $Cl_{n},$ can be generated from $\mathbb{R}^n$ by considering the
relationship $${x}^{2}=-\|{x}\|^{2}$$ for each
${x}\in \mathbb{R}^n$.  We have $\mathbb{R}^n\subseteq Cl_{n}$. If $e_1,\ldots, e_n$ is an orthonormal basis for $\mathbb{R}^n$, then ${x}^{2}=-\|{x}\|^{2}$ tells us that $e_i e_j + e_j e_i= -2\delta_{ij},$ where $\delta_{ij}$ is the Kronecker delta function. Let $A=\{j_1, \cdots, j_r\}\subset \{1, 2, \cdots, n\}$ and $1\leq j_1< j_2 < \cdots < j_r \leq n$. An arbitrary element of the basis of the Clifford algebra can be written as $e_A=e_{j_1}\cdots e_{j_r}.$ Hence for any element $a\in Cl_{n}$, we have $a=\sum_Aa_Ae_A,$ where $a_A\in \mathbb{R}.$

We define the Clifford conjugation as the following:
 $$\bar{a}=\sum_A(-1)^{|A|(|A|+1)/2}a_Ae_A$$
satisfying $\overline{e_{j_1}\cdots e_{j_r}}=(-1)^r e_{j_r}\cdots e_{j_1}$ and $\overline{ab}= \bar{b}\bar{a}$ for $a, b \in Cl_n.$

For each $a=a_0+\cdots +a_{1\cdots n}e_1\cdots e_n\in Cl_n$ the scalar part of $\bar{a}a$ gives the square of the norm of $a,$ namely $a_0^2+\cdots +a_{1\cdots n}^2$\,.

The reversion is given by
 $$\tilde{a}=\sum_A(-1)^{|A|(|A|-1)/2}a_Ae_A,$$ where $|A|$ is the cardinality of $A$. In particular, $\widetilde{e_{j_1}\cdots e_{j_r}}=e_{j_r}\cdots e_{j_1}.$ Also $\widetilde{ab}=\tilde{b}\tilde{a}$ for $a, b \in Cl_{n+1}.$

The Pin and Spin groups play an important role in Clifford analysis.  The Pin group can be defined as
 $$Pin(n): =\{a\in Cl_n : a=y_1 \ldots y_p:
{y_1,\ldots , y_p}\in \mathbb{S}^{n-1}, p\in \mathbb{N}\}$$
and is clearly a group under multiplication in $Cl_n$.

 Now suppose that $y\in \mathbb{S}^{n-1}\subseteq \mathbb{R}^n$. Look at $yxy=yx^{\parallel _y}y+yx^{\perp_y}y=-x^{\parallel _y}+x^{\perp_y}$ where $x^{\parallel _y}$ is the projection of $x$ onto $y$
and $x^{\perp_y}$ is perpendicular to $y$. So $yxy$ gives a reflection of $x$ in the $y$ direction. By the Cartan$-$Dieudonn\'{e} Theorem each $O \in O(n)$  is the composition of a finite number of reflections. If  $a=y_1\ldots y_p\in Pin(n)$, then $\tilde{a}:=y_p\ldots y_1$ and $ax\tilde{a}=O_a(x)$ for some $O_a\in O(n).$ Choosing $y_1, \ldots, y_p$ arbitrarily in $\mathbb{S}^{n-1}$,  we see that the group homomorphism $$\theta: Pin(n)\longrightarrow O(n): a\longmapsto O_a$$ with $a=y_1\ldots y_p$
and $O_a(x)=ax\tilde{a}$ is surjective. Further $-ax(-\tilde{a})=ax\tilde{a}$, so $1, -1\in ker(\theta)$. In fact $ker(\theta)=\{\pm 1\}.$
The Spin group is defined as
$$
Spin(n):=\{a\in Pin(n): a=y_1\ldots y_ p \mbox{ and } p \mbox{ even}\}
$$
and is a subgroup of $Pin(n)$. There is a group homomorphism
 $$\theta: Spin(n)\longrightarrow SO(n)$$ which is surjective with kernel $\{1, -1\}$. See \cite{P} for details.

The Dirac Operator in $\mathbb{R}^n$ is defined to be $$D :=\sum_{j=1}^{n} e_j \frac{\partial}{\partial x_j}.$$ Note $D^2=-\Delta_{n},$ where $\Delta_n$ is the Laplacian in $\mathbb{R}^n$.

Let $\mathcal{H}$$_{k}$ be the space of harmonic polynomials homogeneous of degree $k.$ Let $\mathcal{M}_k$ denote the space of $Cl_n-$ valued polynomials, homogeneous of degree $k$ and such that if $p_k\in$ $\mathcal{M}_k$ then $Dp_k=0.$ Such a polynomial is called a left monogenic polynomial homogeneous of degree $k$. Note if $h_k\in$ $\mathcal{H}_k,$ the space of $Cl_n-$ valued harmonic polynomials homogeneous of degree $k$, then $Dh_k\in$ $\mathcal{M}$$_{k-1}$. But $Dup_{k-1}(u)=(-n-2k+2)p_{k-1}(u),$
so $$\mathcal{H}_k=\mathcal{M}_k\bigoplus u\mathcal{M}_{k-1}, h_k=p_k+up_{k-1}.$$
This is the so-called Almansi-Fischer decomposition of $\mathcal{H}$$_k$, where $\mathcal{M}_{k-1}$ is the space of $Cl_n-$ valued left monogenic polynomials, homogeneous of degree $k-1$. See \cite{BDS, R}.
\par Note that if $Dg(u)=0$ then $\bar{g}(u)\bar{D}=-\bar{g}(u)D=0$. So we can talk of right monogenic polynomials, homogeneous of degree $k$ and we obtain by conjugation a right Almansi-Fisher decomposition, $$\mathcal{H}_k=\overline{\mathcal{M}_k}\bigoplus \overline{\mathcal{M}}_{k-1}u,$$ where $\overline{\mathcal{M}_k}$ stands for the space of right monogenic polynomials homogeneous of degree $k$.

Let $P_k$ be the left projection map
 $$P_k:  \mathcal{H}_k\rightarrow \mathcal{M}_k,$$ then the left Rarita-Schwinger operator $R_k$ is defined by (see \cite{BSSV,BSSV1,DLRV,Va1,Va2})
$$R_kg(x,u)=P_kD_xg(x,u),$$
where $D_x$ is the Dirac operator with respect to $x$ and $g(x,u): U\times \mathbb{R}^n\to Cl_n$ is a monogenic polynomial homogeneous of degree $k$ in $u,$ and $U$ is a domain in $\mathbb{R}^n$.
The left Rarita-Schwinger equation is defined to be
$$
R_k g(x,u)=0.
$$
We also have a right projection $P_{k,r}: \mathcal{H}_k \rightarrow \overline{\mathcal{M}_k},$ and a right Rarita-Schwinger equation $g(x,u)D_xP_{k,r}=g(x,u)R_k=0.$

A M\"{o}bius transformation is a finite composition of orthogonal transformations, inversions, dilations, and translations.
Ahlfors \cite{A} and Vahlen \cite{V} show that given a M\"{o}bius transformation $y=\phi(x)$ on $\mathbb{R}^n\bigcup\{\infty\}$ it can be expressed as $y=(ax+b)(cx+d)^{-1}$ where $a, b,c,d\in Cl_n$ and satisfy the following conditions:
\begin{enumerate}
 \item $a, b, c, d$ are all products of vectors in $\mathbb{R}^n.$
 \item $a\tilde b, c\tilde d, \tilde bc, \tilde da \in\mathbb{R}^n.$
 \item $a\tilde d-b\tilde c=\pm1.$
\end{enumerate}
When $c=0, \phi(x)=(ax+b)(cx+d)^{-1}=axd^{-1}+bd^{-1}=\pm ax\tilde a+bd^{-1}.$ Now assume $c\neq 0,$ then $\phi(x)=(ax+b)(cx+d)^{-1}=ac^{-1}\pm(cx\tilde{c}+d\tilde{c})^{-1}.$ These are the so-called Iwasawa decompositions. Using this notation and the conformal weights, $f(\phi(x))$ is changed to $J(\phi,x)f(\phi(x)),$ where $J(\phi,x)=\di\frac{\widetilde{cx+d}}{\|cx+d\|^n}$. Note when $\phi(x)=x+a$ then $J(\phi,x)\equiv 1.$

\section{The $Q_k$ operators and their kernels}

As $$I-P_k:  \mathcal{H}_k\rightarrow u\mathcal{M}_{k-1},$$ where $I$ is the identity map, then we can define the left remaining operators
$$Q_k:=(I-P_k)D_x:u\mathcal{M}_{k-1}\to u\mathcal{M}_{k-1}\quad uf(x,u):\to (I-P_k)D_xuf(x,u).$$ See\cite{BSSV}.

The left remaining equation is defined to be $(I-P_k)D_xuf(x,u)=0$ or $Q_kuf(x,u)=0,$ for each $x$ and $(x,u)\in U\times \mathbb{R}^n$, where $U$ is a domain in $\mathbb{R}^n$ and $f(x,u)\in \mathcal{M}_{k-1}.$

We also have a right remaining operator $$Q_{k,r}:=D_x(I-P_{k,r}):\overline{\mathcal{M}}_{k-1}u\to \overline{\mathcal{M}}_{k-1}u\quad g(x,u)u:\to g(x,u)uD_x(I-P_{k,r}),$$ where $g(x,u)\in \overline{\mathcal{M}}_{k-1}$.

Consequently, the right remaining equation is $g(x,u)uD_x(I-P_{k,r})=0$ or $g(x,u)uQ_{k,r}=0.$

Now let us establish the conformal invariance of the remaining equation $Q_kuf(x,u)=0$.

It is easy to see that $I-P_k$ is conformally invariant under the M\"{o}bius transformations, since the projection operator $P_k$ is conformally invariant (see \cite{BSSV, DLRV}). By considering orthogonal transformation, inversion, dilation and translation and applying the same arguments in \cite{DLRV} used to establish the intertwining operators for Rarita-Schwinger operators, we can easily obtain the intertwining operators for $Q_k$ operators:
\begin{theorem}
$$J_{-1}(\phi,x)Q_{k,u}uf(y,u)=Q_{k,w}wJ(\phi,x)f(\phi(x),\di\frac{\widetilde{(cx+d)}w(cx+d)}{\|cx+d\|^2}),$$
where $Q_{k,u}$ and $Q_{k,w}$ are the remaining operators with respect to $u$ and $w$ respectively, $y=\phi(x)$ is the M\"{o}bius transformation, $J(\phi,x)=\di\frac{\widetilde{cx+d}}{\|cx+d\|^n},J_{-1}(\phi,x)=\di\frac{cx+d}{\|cx+d\|^{n+2}},$ and $u=\di\frac{\widetilde{(cx+d)}w(cx+d)}{\|cx+d\|^2}$ for some $w\in \mathbb{R}^n.$
\end{theorem}

Consequently, we have \\
$Q_{k,u}uf(x,u)=0$ implies $Q_{k,w}wJ(\phi,x)f(\phi(x),\di\frac{\widetilde{(cx+d)}w(cx+d)}{\|cx+d\|^2})=0$. This tells us that the remaining equation $Q_kuf(x,u)=0$ is conformally invariant under M\"{o}bius transformations.

The reproducing kernel of $\mathcal{M}_k$ with respect to integration over $\mathbb{S}^{n-1}$ is given by (see \cite{BDS,DLRV}) $$Z_k(u,v):=\di\sum_\sigma P_\sigma(u)V_\sigma(v)v,$$ where
$$
P_\sigma(u)=\di\frac{1}{k!}\di\Sigma(u_{i_1}-u_1e_1^{-1}e_{i_1})\ldots(u_{i_k}-u_1e_1^{-1}e_{i_k}),
V_\sigma(v)=\di\frac{\partial^kG(v)}{\partial v_{2}^{j_2}\ldots \partial v_{n}^{j_{n}}}\,
$$
$j_2+\ldots +j_{n}=k,~~\mbox{and}~~ i_k\in\{2,\cdots,n\}.$ Here summation is taken over all permutations of the monomials without repetition.  This function is left monogenic in $u$ and right monogenic polynomial in $v$ and it is homogeneous of degree $k$. See \cite{BDS} and elsewhere.

Let us consider the polynomial $uZ_{k-1}(u,v)v$ which is harmonic, homogeneous degree of $k$ in both $u$ and $v$. Since $uZ_{k-1}(u,v)v$ does not depend on $x$, $Q_kuZ_{k-1}(u,v)v=0$.

Now applying inversion from the left, we obtain $$H_k(x,u,v):=\di\frac{-1}{\omega_n c_k}u\di\frac{x}{\|x\|^n}Z_{k-1}(\di\frac{xux}{\|x\|^2},v)v$$ is a non-trivial  solution to $Q_kuf(x,u)=0,$ where $c_k=\di\frac{n-2}{n-2+2k}.$

 Similarly, applying inversion from the right, we obtain $$\di\frac{-1}{\omega_n c_k}uZ_{k-1}(u,\di\frac{xvx}{\|x\|^2})\di\frac{x}{\|x\|^n}v$$ is a non-trivial solution to $f(x,v)vQ_{k,r}=0.$ Using the similar arguments in \cite{DLRV}, we can show that two representations of the solutions are equal. The details are given in the following
 $$\begin{array}{ll}
 \di\frac{-1}{\omega_n c_k}uZ_{k-1}(u,\di\frac{xvx}{\|x\|^2})\di\frac{x}{\|x\|^n}v\\
 \\
 =\di\frac{-1}{\omega_n c_k}u\di\frac{-x}{\|x\|}Z_{k-1}(\di\frac{xux}{\|x\|^2},v)\di\frac{x}{\|x\|}\di\frac{x}{\|x\|^n}v
 =\di\frac{-1}{\omega_n c_k}u\di\frac{x}{\|x\|^n}Z_{k-1}(\di\frac{xux}{\|x\|^2},v)v.
 \end{array}$$
 In fact $H_k(x,u,v)$ is the fundamental solution to the $Q_k$ operator.

\section{Some basic integral formulas related to $Q_k$ operators}
In this section, we will establish some basic integral formulas associated with $Q_k$ operators.
\begin{definition}\cite{DLRV} \quad For any $Cl_n-$valued polynomials $P(u), Q(u)$, the inner product $(P(u), Q(u))_u$ with respect to $u$ is given by $$(P(u), Q(u))_u=\di\int_{\mathbb{S}^{n-1}}P(u)Q(u)ds(u).$$
\end{definition}
For any $p_k \in \mathcal{M}_k,$ one obtains (see \cite{BDS})
$$
p_k(u)=(Z_k(u,v), p_k(v))_v=\int_{\mathbb{S}^{n-1}}Z_k(u,v)p_k(v)ds(v).
$$

\par Now if we combine Stokes' Theorems of the Dirac operator and the Rarita-Schwinger operator, then we have two versions of Stokes' Theorem for the $Q_k$ operators .
\begin{theorem}(Stokes' Theorem for $Q_k$ operators) Let $\Omega'$ and $\Omega$ be domains in $\mathbb{R}^n$ and suppose the closure of $\Omega$ lies in $\Omega'$. Further suppose the closure of $\Omega$ is compact and the boundary of $\Omega,$ $\partial\Omega$, is piecewise smooth. Then for $f, g \in C^1(\Omega',$$\mathcal{M}_k)$, we have version 1
$$\begin{array}{ll}
\di\int_\Omega[(g(x,u)Q_{k,r}, f(x,u))_u+(g(x,u), Q_kf(x,u))_u]dx^n\\
\\
=\di\int_{\partial\Omega}\left(g(x,u), (I-P_k)d\sigma_xf(x,u)\right)_u\\
\\
=\di\int_{\partial\Omega}\left(g(x,u)d\sigma_x(I-P_{k,r}), f(x,u)\right)_u.
\end{array}$$
Then for $f, g \in C^1(\Omega',$$\mathcal{M}_{k-1})$, we have version 2
$$\begin{array}{ll}
\di\int_\Omega[(g(x,u)uQ_{k,r}), uf(x,u))_u+(g(x,u)u, Q_kuf(x,u))_u]dx^n\\
\\
=\di\int_{\partial\Omega}\left(g(x,u)u, (I-P_k)d\sigma_xuf(x,u)\right)_u\\
\\
=\di\int_{\partial\Omega}\left(g(x,u)ud\sigma_x(I-P_{k,r}), uf(x,u)\right)_u.
\end{array}$$

\end{theorem}
\par {\bf Proof:}\quad It is easy to get version 1 of Stokes' Theorem for the $Q_k$ operators by combining Stokes' Theorems of the Dirac operator and the Rarita-Schwinger operators.

Now we shall prove version 2 of Stokes' Theorem.

First of all, we want to prove that $$\di\int_{\partial\Omega}\left(g(x,u)u, (I-P_k)d\sigma_xuf(x,u)\right)_u
=\di\int_{\partial\Omega}\left(g(x,u)ud\sigma_x(I-P_{k,r}), uf(x,u)\right)_u.$$
Here $d\sigma_x=n(x)d\sigma(x)$.
By the Almansi-Fischer decomposition, we have $$g(x,u)un(x)uf(x,u)=g(x,u)u[f_1(x,u)+uf_2(x,u)]=[g_1(x,u)+g_2(x,u)u]uf(x,u),$$
so $$\begin{array}{ll}g(x,u)ud\sigma_xuf(x,u)\\
\\
=g(x,u)u[f_1(x,u)+uf_2(x,u)]d\sigma(x)=[g_1(x,u)+g_2(x,u)u]uf(x,u)d\sigma(x),\end{array}$$where $f_1(x,u), f_2(x,u), g_1(x,u), g_2(x,u)$ are left or right monogenic polynomials in $u.$ Now integrating the above formula over the unit sphere in $\mathbb{R}^n$, one gets
$$ \begin{array}{ll}\di\int_{\mathbb{S}^{n-1}}g(x,u)ud\sigma_xuf(x,u)ds(u)\\
\\
=\di\int_{\mathbb{S}^{n-1}}g(x,u)uuf_2(x,u)d\sigma(x)ds(u)= \di\int_{\mathbb{S}^{n-1}}g_2(x,u)uuf(x,u)d\sigma(x)ds(u).\end{array}$$ This follows from the fact that $$\di\int_{\mathbb{S}^{n-1}}g(x,u)uf_1(x,u)ds(u)=\di\int_{\mathbb{S}^{n-1}}g_1(x,u)uf(x,u)ds(u)=0.$$ See \cite{BDS}.

Thus \begin{eqnarray*}\di\int_{\partial\Omega}\left(g(x,u)u, (I-P_k)d\sigma_xuf(x,u)\right)_u&=&\di\int_{\partial\Omega} \int_{\mathbb{S}^{n-1}}g(x,u)u((I-P_k)d\sigma_xuf(x,u))ds(u)\\
&=&\di\int_{\partial\Omega} \int_{\mathbb{S}^{n-1}}g(x,u)u uf_2(x,u)ds(u)d\sigma(x)\\
&= &\di\int_{\partial\Omega}\int_{\mathbb{S}^{n-1}} g_2(x,u)u u f(x,u)ds(u)d\sigma(x)\\
&=&\di\int_{\partial\Omega} \int_{\mathbb{S}^{n-1}}(g(x,u)d\sigma_x(I-P_{k,r}))uf(x,u)ds(u)\\
&=&\di\int_{\partial\Omega}\left(g(x,u)ud\sigma_x(I-P_{k,r}), uf(x,u)\right)_u.
\end{eqnarray*}
Secondly, we need to show $$\begin{array}{ll}\di\int_\Omega[(g(x,u)uQ_{k,r}, uf(x,u))_u+(g(x,u)u, Q_kuf(x,u))_u]dx^n\\
\\
=\di\int_{\partial\Omega}\left(g(x,u)u, (I-P_k)d\sigma_xuf(x,u)\right)_u.\end{array}$$
Consider the integral
\begin{eqnarray}\label{1}\di\int_\Omega[(g(x,u)uD_xP_{k,r},uf(x,u))_u+(g(x,u)u,P_kD_xuf(x,u))_u]dx^n\nonumber\\
 \nonumber\\
=\di\int_\Omega\int_{\mathbb{S}^{n-1}}[(g(x,u)uD_xP_{k,r})uf(x,u)+g(x,u)u(P_kD_xuf(x,u))]ds(u)dx^n. \end{eqnarray}

Since $g(x,u)uD_xP_{k,r}, f(x,u), g(x,u)$ and $P_kD_xuf(x,u)$ are monogenic functions in $u$,
$$\int_{\mathbb{S}^{n-1}}(g(x,u)uD_xP_{k,r})uf(x,u)ds(u)=0=\int_{\mathbb{S}^{n-1}}g(x,u)u(P_kD_xuf(x,u))ds(u).$$
Thus the previous integral (\ref{1}) equals zero.

By Stokes' Theorem for the Dirac operator, we have
$$\begin{array}{lll}\di\int_\Omega[(g(x,u)uD_x,uf(x,u))_u+(g(x,u)u,D_xuf(x,u))_u]dx^n\\
\\
=\di\int_\Omega\int_{\mathbb{S}^{n-1}}[(g(x,u)uD_x)uf(x,u)+g(x,u)u(D_xuf(x,u))]ds(u)dx^n\\
\\
=\di\int_{\partial\Omega}\int_{\mathbb{S}^{n-1}}[(g(x,u)ud\sigma_xuf(x,u))]ds(u)\\
\\
=\di\int_{\partial\Omega}(g(x,u)u,d\sigma_xuf(x,u))_u.\end{array}$$
But $$\di\int_{\partial\Omega}(g(x,u)u,P_kd\sigma_xuf(x,u))_u=\di\int_{\partial\Omega}\int_{\mathbb{S}^{n-1}}g(x,u)u (P_kd\sigma_xuf(x,u))ds(u)=0,$$
since $\di\int_{\mathbb{S}^{n-1}}g(x,u)u (P_kd\sigma_xuf(x,u))ds(u)=0.$

Therefore we have shown

$$\begin{array}{ll}\di\int_\Omega[(g(x,u)uQ_{k,r}, uf(x,u))_u+(g(x,u)u, Q_kuf(x,u))_u]dx^n\\
\\
=\di\int_{\partial\Omega}\left(g(x,u)u, (I-P_k)d\sigma_xuf(x,u)\right)_u. \quad \blacksquare\end{array}$$

\begin{remark}
In the  proof of the previous theorem it is proved that
\begin{eqnarray}\label{2}
\di\int_{\partial\Omega}\left(g(x,u)u, (I-P_k)d\sigma_xuf(x,u)\right)_u=\di\int_{\partial\Omega}\left(g(x,u)u, d\sigma_xuf(x,u)\right)_u.
\end{eqnarray}
\end{remark}

\begin{theorem}(Borel-Pompeiu Theorem)Let $\Omega'$ and $\Omega$ be as in the previous Theorem. Then for $f\in C^1(\Omega',$$\mathcal{M}_{k-1})$ and $y\in \Omega,$ we obtain
$$\begin{array}{ll}
uf(y,u)=\di\int_\Omega(H_k(x-y,u,v),Q_kvf(x,v))_vdx^n\\
\\
-\di\int_{\partial\Omega}\left(H_k(x-y,u,v), (I-P_k)d\sigma_xvf(x,v)\right)_v.
\end{array}$$ Here we will use the representation $H_k(x-y,u,v)=\di\frac{-1}{\omega_n c_k}uZ_{k-1}(u,\di\frac{(x-y)v(x-y)}{\|x-y\|^2})\di\frac{x-y}{\|x-y\|^n}v$. \end{theorem}

\par{\bf Proof:}\qquad  Consider a ball $B(y,r)$ centered at $y$ with radius $r$ such that $\overline{B(y,r)}\subset \Omega$. We have
$$\begin{array}{ll}
\di\int_\Omega(H_k(x-y,u,v),Q_kvf(x,v))_vdx^n\\
\\
=\di\int_{\Omega\setminus{B(y,r)}}(H_k(x-y,u,v),Q_kvf(x,v))_vdx^n\\
\\
+\di\int_{B(y,r)}(H_k(x-y,u,v),Q_kvf(x,v))_vdx^n.
\end{array}$$
\par The last integral in the previous equation tends to zero as $r$ tends to zero. This follows from the degree of homogeneity of $x-y$ in $H_k(x-y,u,v)$. Now applying Stokes' Theorem version 2 to the first integral, one gets
$$\begin{array}{ll}\di\int_{\Omega\setminus{B(y,r)}}(H_k(x-y,u,v),Q_kvf(x,v))_vdx^n\\
\\
=\di\int_{\partial\Omega}(H_k(x-y,u,v),(I-P_k)d\sigma_xvf(x,v))_v-\di\int_{\partial {B(y,r)}}(H_k(x-y,u,v),(I-P_k)d\sigma_xvf(x,v))_v.
 \end{array}$$
\par Now let us look at the integral $$\begin{array}{ll}\di\int_{\partial {B(y,r)}}(H_k(x-y,u,v),(I-P_k)d\sigma_xvf(x,v))_vdx^n\\
\\
=\di\int_{\partial {B(y,r)}}(H_k(x-y,u,v),(I-P_k)d\sigma_xvf(y,v))_v\\
\\
+\di\int_{\partial {B(y,r)}}(H_k(x-y,u,v),(I-P_k)d\sigma_xv[f(x,v)-f(y,v)])_v.
 \end{array}$$
 Since the second integral on the right hand side tends to zero as $r$ goes to zero because of the continuity of $f$, we only need to deal with the first integral
 $$\begin{array}{llll}\di\int_{\partial {B(y,r)}}(H_k(x-y,u,v),(I-P_k)d\sigma_xvf(y,v))_v \\
 \\
 =\di\int_{\partial{B(y,r)}}\int_{\mathbb{S}^{n-1}}H_k(x-y,u,v)(I-P_k)d\sigma_xvf(y,v)ds(v)\\
 \\
 =\di\int_{\partial{B(y,r)}}\int_{\mathbb{S}^{n-1}}\di\frac{-1}{\omega_n c_k}uZ_{k-1}\left(u,\di\frac{(x-y)v(x-y)}{\|x-y\|^2}\right)\di\frac{x-y}{\|x-y\|^n}v(I-P_k)n(x)vf(y,v)ds(v)d\sigma(x),
 \end{array}$$ where $n(x)$ is the unit outer normal vector and $d\sigma(x)$ is the scalar measure on $\partial B(y,r).$ Now $n(x)$ here is $\di\frac{y-x}{\|x-y\|}.$ Using equation (\ref{2}) the previous integral becomes
 $$\begin{array}{ll}
 \di\int_{\partial{B(y,r)}}\int_{\mathbb{S}^{n-1}}\di\frac{1}{\omega_n c_k}uZ_{k-1}\left(u,\di\frac{(x-y)v(x-y)}{\|x-y\|^2}\right)\di\frac{x-y}{\|x-y\|^n}v\di\frac{x-y}{\|x-y\|}
vf(y,v) ds(v)d\sigma(x)\\
\\
=\di\frac{1}{\omega_n c_k}\di\int_{\partial{B(y,r)}}\di\frac{1}{r^{n-1}}\int_{\mathbb{S}^{n-1}}uZ_{k-1}\left(u,\di\frac{(x-y)v(x-y)}{\|x-y\|^2}\right)\di\frac{x-y}{\|x-y\|}v\di\frac{x-y}{\|x-y\|}
vf(y,v) ds(v)d\sigma(x)
\end{array}$$
 \par Since $Z_{k-1}\left(u,\di\frac{(x-y)v(x-y)}{\|x-y\|^2}\right)\di\frac{x-y}{\|x-y\|}v\di\frac{x-y}{\|x-y\|}$ is a harmonic polynomial with degree $k$ in $v$, we can apply Lemma $5$ in \cite{DLRV}, then the integral is equal to
$$\begin{array}{ll}
\di\int_{\mathbb{S}^{n-1}}uZ_{k-1}(u,v)vvf(y,v)ds(v)d\sigma(x)\\
\\
=-u\di\int_{\mathbb{S}^{n-1}}Z_{k-1}(u,v)f(y,v)ds(v)=-uf(y,u),
\end{array}$$
\par Therefore, when the radius $r$ tends to zero, we obtain the desired result. \qquad $\blacksquare$

\par Now if the function has compact support in $\Omega$, then by Borel-Pompeiu Theorem we obtain:\\
\begin{theorem} $\di\iint_{\mathbb{R}^n}(H_k(x-y,u,v),Q_kv\phi(x,v))_vdx^n=u\phi(y,u)$ for each $\phi\in C_0^{\infty}(\mathbb{R}^n)$.\end{theorem}

\par Now suppose $vf(x,v)$ is a solution to the $Q_k$ operator, then using the Borel-Pompeiu Theorem we have:\\

\begin{theorem} (Cauchy Integral Formula) If $Q_kvf(x,v)=0,$ then for $y\in \Omega,$
$$\begin{array}{ll}
uf(y,u)=-\di\int_{\partial\Omega}\left(H_k(x-y,u,v), (I-P_k)d\sigma_xvf(x,v)\right)_v\\

=-\di\int_{\partial\Omega}\left(H_k(x-y,u,v)d\sigma_x(I-P_{k,r}), vf(x,v)\right)_v.
\qquad \blacksquare \end{array}$$ \end{theorem}

 We also can talk about a Cauchy transform for the $Q_k$ operators:

\begin{definition}\quad For a domain $\Omega\subset \mathbb{R}^n$ and a function $f: \Omega\times\mathbb{R}^n \longrightarrow Cl_n,$ where $f(x,u)$ is monogenic in $u$ with degree $k-1$, the Cauchy $($or $T_k$-transform$)$ of $f$ is formally defined to be $$(T_kvf)(y,v)=-\iint_\Omega \left(H_k(x-y,u,v), uf(x,u)\right)_udx^n,\qquad y\in \Omega.$$\end{definition}

\begin{theorem} \quad   $$Q_k\di\iint_{\mathbb{R}^n}\left(H_k(x-y,u,v), v\phi(x,v)\right)_vdx^n=u\phi(y,u), \mbox{for}~~ \phi\in C_0^{\infty}(\mathbb{R}^n).$$ Here we use the representation $H_k(x-y,u,v)=\di\frac{-1}{\omega_n c_k}uZ_{k-1}(u,\di\frac{(x-y)v(x-y)}{\|x-y\|^2})\di\frac{x-y}{\|x-y\|^n}v.$ \end{theorem}
\par {\bf Proof:} \quad  For each fixed $y\in \mathbb{R}^n$, we can construct a bounded rectangle $R(y)$ centered at $y$ in $\mathbb{R}^n$.
\par Then

$$\begin{array}{lll}Q_k\di\iint_{\mathbb{R}^n\setminus R(y)}(H_k(x-y,u,v), v\phi(x,v))_vdx^n\\
\\
=(I-P_k)D_y\di\iint_{\mathbb{R}^n\setminus R(y)}(H_k(x-y,u,v), v\phi(x,v))_vdx^n=0.\end{array}$$

 Now consider $$\begin{array}{lll}
\di\frac{\partial}{\partial y_i}\di\iint_{R(y)}(H_k(x-y,u,v), v\phi(x,v))_vdx^n\\
\\
=\lim_{\varepsilon\to 0}\di\frac{1}{\varepsilon}[\di\iint_{R(y)}(H_k(x-y,u,v), v\phi(x,v))_vdx^n\\
\\
-\di\iint_{R(y+\varepsilon e_i)}(H_k(x-y-\varepsilon e_i,u,v), v\phi(x,v))_vdx^n]\end{array}$$
$$\begin{array}{lll}=\di\iint_{R(y)}(H_k(x-y,u,v), \di\frac{\partial v\phi(x,v)}{\partial x_i})_vdx^n\\
\\
+\di\int_{\partial R_1(y)\cup \partial R_2(y)}(H_k(x-y,u,v), v\phi(x,v))_vd\sigma(x)
\end{array}$$
where $\partial R_1(y)$ and $\partial R_2(y)$ are the two faces of $R(y)$ with normal vectors $\pm e_i.$
So $$\begin{array}{ll}D_y\di\iint_{R(y)}(H_k(x-y,u,v), v\phi(x,v))_vdx^n\\
\\
=\di\iint_{R(y)}\sum_{i=1}^n e_i(H_k(x-y,u,v), \di\frac{\partial {v\phi(x,v)}}{\partial {x_i}})_vdx^n\\
\\
+\di\int_{\partial R(y)}n(x)(H_k(x-y,u,v),v\phi(x,v))_vd\sigma(x).\end{array}$$

The first integral tends to zero as the volume of $R(y)$ tends to zero. Thus we will pay attention to the integral
$$(I-P_k)\di\int_{\partial R(y)}n(x)(H_k(x-y,u,v),v\phi(x,v))_vd\sigma(x).$$
\par This is equal to
$$(I-P_k)\di\int_{\partial R(y)}\di\int_{\mathbb{S}^{n-1}}n(x)H_k(x-y,u,v)v\phi(x,v)ds(v)d\sigma(x),$$ which in turn is equal to
$$\begin{array}{ll}
(I-P_k)\di\int_{\partial R(y)}\di\int_{\mathbb{S}^{n-1}}n(x)H_k(x-y,u,v)v\phi(y,v)ds(v)d\sigma(x)\\
\\
+(I-P_k)\di\int_{\partial R(y)}\di\int_{\mathbb{S}^{n-1}}n(x)H_k(x-y,u,v)v(\phi(x,v)-\phi(y,v))ds(v)d\sigma(x).\end{array}$$
But the last integral on the right side of the above formula tends to zero as the surface area of $\partial R(y)$ tends to zero because of the degree of the homogeneity of $x-y$ in $H_k$ and the continuity of the function $\phi$. Hence we are left with
$$(I-P_k)\di\int_{\partial R(y)}\di\int_{\mathbb{S}^{n-1}}n(x)H_k(x-y,u,v)v\phi(y,v)ds(v)d\sigma(x).$$
\par By Stokes' Theorem this is equal to $$(I-P_k)\di\int_{\partial B(y,r)}\di\int_{\mathbb{S}^{n-1}}n(x)H_k(x-y,u,v)v\phi(y,v)ds(v)d\sigma(x).$$  In turn this is equal to
$$\begin{array}{ll}
(I-P_k)\di\int_{\partial B(y,r)}\di\int_{\mathbb{S}^{n-1}}-\di\frac{x-y}{\|x-y\|}\di\frac{-1}{\omega_n c_k}uZ_{k-1}(u,\di\frac{(x-y)v(x-y)}{\|x-y\|^2})\di\frac{x-y}{\|x-y\|^n}vv\phi(y,v)ds(v)d\sigma(x)\\
\\
=(I-P_k)\di\int_{\partial B(y,r)}\di\frac{-1}{\omega_n c_k}\di\int_{\mathbb{S}^{n-1}}\di\frac{x-y}{\|x-y\|}uZ_{k-1}(u,\di\frac{(x-y)v(x-y)}{\|x-y\|^2},v)\di\frac{x-y}{\|x-y\|^n}\phi(y,v)ds(v)d\sigma(x).
\end{array}$$
Since $Z_{k-1}(u,v)$ is the reproducing kernel of $\mathcal{M}_{k-1}$, $\pm\tilde{a}Z_{k-1}(au\tilde{a},av\tilde{a})a$ is also the reproducing kernel of $\mathcal{M}_{k-1}$ for each $a\in Pin(n)$. See \cite{DLRV}. Now let $a=\di\frac{x-y}{\|x-y\|},$ the previous integral equals
$$\begin{array}{ll}
(I-P_k)\di\int_{\partial B(y,r)}\di\frac{1}{\omega_n c_k}\di\int_{\mathbb{S}^{n-1}}\di\frac{x-y}{\|x-y\|}u\di\frac{x-y}{\|x-y\|}Z_{k-1}(\di\frac{(x-y)u(x-y)}{\|x-y\|^2},v)\di\frac{x-y}{\|x-y\|}\di\frac{x-y}{\|x-y\|^n}\\
\\
\phi(y,v)ds(v)d\sigma(x)\\
\\
=(I-P_k)\di\int_{\partial B(y,r)}\di\int_{\mathbb{S}^{n-1}}\di\frac{1}{\omega_n c_k}\di\frac{1}{r^{n-1}}\di\frac{x-y}{\|x-y\|}u\di\frac{x-y}{\|x-y\|}Z_{k-1}(\di\frac{(x-y)u(x-y)}{\|x-y\|^2},v)\phi(y,v)ds(v)d\sigma(x).
\end{array}$$
Applying Lemma 5 in \cite{DLRV}, the integral becomes
$$
(I-P_k)\di\int_{\mathbb{S}^{n-1}}uZ_{k-1}(u,v)\phi(y,v)ds(v)=(I-P_k)u\phi(y,u)=u\phi(y,u).\quad \blacksquare
$$

%%%%%%%%%%%%%%%%%%%%%%%%%%
 \section{The $Q_k$ operators on the sphere}
In this section, we will extend the results for the $Q_k$ operators in $\mathbb{R}^n$ from the previous sections to the sphere.

 Consider the Cayley transformation
$ C : \mathbb{R}^n \to \mathbb{S}^n $, where $\mathbb{S}^n$ is the unit sphere in $\mathbb{R}^{n+1}$,
defined by $C(x)= (e_{n+1}x +1)(x + e_{n+1})^{-1} $, where $ x= x_1e_1 + \cdots + x_ne_n \in \mathbb{R}^n $,
and $e_{n+1}$ is a unit vector in $\mathbb{R}^{n+1}$ which is orthogonal to $\mathbb{R}^n $.
Now $C(\mathbb{R}^n) = \mathbb{S}^n \setminus\{e_{n+1}\}$. Suppose $ x_s \in \mathbb{S}^n $ and
$x_s = x_{s_1}e_1 + \cdots + x_{s_n}e_n + x_{s_{n+1}}e_{n+1}$, then we have
$C^{-1}(x_s) = (-e_{n+1}x_s +1)(x_s - e_{n+1})^{-1}$.
\\

The Dirac operator over the $n$-sphere $\mathbb{S}^n$ has the form $D_{s} = w(\Gamma + \frac{n}{2})$, where
$w \in \mathbb{S}^n $ and $ \Gamma =
\displaystyle\sum_{i<j, i=1}^{n}{e_ie_j(w_i\frac{\partial}{\partial w_j} - w_j\frac{\partial}{\partial w_i})}$. See \cite{CM,LR,R1,R2,Va3}.

Let $U$ be a domain in $\mathbb{R}^n$. Consider a function $f_{\star}: U \times \mathbb{R}^{n} \to Cl_{n+1}$ such that for each $x\in U$,
$f_{\star}(x, u)$ is a left monogenic polynomial homogeneous of degree $k-1$ in $u$.
This function reduces to $f(x_s, u)$ on $C(U) \times \mathbb{R}^{n} $ and $f(x_s, u)$ takes its values in $Cl_{n+1},$
where $C(U)\subset \mathbb{S}^n$ and $f(x_s, u)$ is a left monogenic polynomial homogeneous of degree $k-1$ in $u$.

The
 left $n$-spherical remaining operator on the sphere is defined to be
$$
Q_k^S=:(I-P_k)D_{s,x_s},
$$ where $D_{s,x_s}$ is the Dirac operator on the sphere with respect to $x_s$.

Hence the left $n$-spherical remaining equation is defined to be $$Q_k^Suf(x_s, u)=0.$$
On the other hand, the  right $n$-spherical remaining operator is defined to be
$$
Q_{k,r}^S:=D_{s,x_s}(I-P_{k,r}).
$$
The right $n$-spherical remaining equation is defined to be $g(x_s,v)vQ_{k,r}^S=0,$ where $g(x_s,v)\in \overline{\mathcal{M}}_{k -1}.$

\subsection{The intertwining operators for $Q_k^S$ and $Q_k$ operators and the conformal invariance of $Q_k^Suf(x_s,u)=0$}

First let us recall that if $ f(u)\in \mathcal{M}_{k -1}$ then it trivially extends to $F(v) = f(u + u_{n+1}e_{n+1})$
with $u_{n+1}\in \mathbb{R}$ and $F(v)= f(u)$ for all $u_{n+1}\in \mathbb{R}$. Consequently $D_{n+1}F(v)= 0$
where $D_{n+1}= \displaystyle\sum_{j=1}^{n+1} {e_j \frac{\partial}{\partial u_j}}$.

If $ f(u)\in \mathcal{M}_{k-1} $ then for any boundary of a piecewise smooth bounded domain
$U \subseteq \mathbb{R}^{n} $ by Cauchy's Theorem

\begin{eqnarray}\label{3}\int_{\partial U}{n(u)f(u)d\sigma(u)} = 0. \end{eqnarray}

Suppose now $a\in \mbox{Pin}(n+1)$ and $u= aw\tilde{a}$ then although $u \in \mathbb{R}^{n}$ in
general $w$ belongs to the hyperplane $a^{-1}\mathbb{R}^{n}\tilde{a}^{-1}$ in $\mathbb{R}^{n+1}$.

By applying a change of variable, up to a sign the integral (\ref{3}) becomes
\begin{eqnarray}\label{4}
\int_{a^{-1}\partial U\tilde{a}^{-1}}{ an(w)\tilde{a}F(aw\tilde{a})d\sigma(w)} = 0.\end{eqnarray}
As $\partial U$ is arbitrary then on applying Stokes' Theorem to (\ref{4}) we see
that
$$
D_a \tilde{a}F(aw\tilde{a}) = 0,~~ \mbox{where} ~~ D_a : = D_{n+1}\bigl|_{a^{-1}\mathbb{R}^{n}\tilde{a}^{-1}}.
$$ See \cite{LRV}.

From now on all functions on spheres take their values in $Cl_{n+1}.$

Now let $f(x_s,u): U_s \times \mathbb{R}^{n}\to Cl_{n+1}$ be a monogenic polynomial
homogeneous of degree $k$ in $u$ for each $x_s \in U_s$, where $U_s$ is a domain in $ \mathbb{S}^n.$

It is known from section 3 that $I-P_k$ is conformally invariant under a general M\"{o}bius transformation
over $\mathbb{R}^{n}$. This trivially extends to M\"{o}bius transformations on $\mathbb{R}^{n+1}$.
It follows that if we restrict $x_s $ to $ \mathbb{S}^n,$ then $I-P_k$ is also conformally invariant under
the Cayley transformation $C$ and its inverse $C^{-1},$ with $x\in \mathbb{R}^n$.

We can use the intertwining formulas for $D_x$ and $D_{s,x_s}$ given in \cite{LR} to establish the intertwining
formulas for $Q_k$ and $Q_k^S.$

\begin{theorem}

$$\begin{array}{ll}
-J_{-1}(C^{-1},x_s)Q_{k,u}uf(x,u)
=Q_{k,w}^SwJ(C^{-1},x_s)f(C^{-1}(x_s),\di\frac{(x_s-e_{n+1})w(x_s-e_{n+1})}{\|x_s-e_{n+1}\|^2}),
\end{array}$$
where $Q_{k,u}$ are the remaining operators with respect to $u \in \mathbb{R}^{n}$, $Q_{k,w}^S$ are the remaining operators on the sphere with respect to $w  \in \mathbb{S}^{n}$,
$u=\displaystyle\frac{(x_s - e_{n+1})w(x_s - e_{n+1})}{||x_s - e_{n+1}||^2},$
$J(C^{-1},x_s)=\di\frac{x_s - e_{n+1}}{\|x_s - e_{n+1}\|^n}$ is the conformal weight for the inverse of the Cayley transformation
and $J_{-1}(C^{-1},x_s)=\di\frac{x_s - e_{n+1}}{\|x_s - e_{n+1}\|^{n+2}}.$
\end{theorem}
{\bf Proof:}\quad  In \cite{LR} it is shown that $D_x=J_{-1}(C^{-1},x_s)^{-1}D_{s,x_s}J(C^{-1},x_s).$

Set $u=\di\frac{J(C^{-1},x_s)wJ(C^{-1},x_s)}{\|J(C^{-1},x_s)\|^2}$ for some $w  \in \mathbb{R}^{n+1}$. Consequently,
$$\begin{array}{ll}Q_{k,u}uf(x,u)=(I-P_{k,u})D_xuf(x,u)\\
\\
=(I-P_{k,u})J_{-1}(C^{-1},x_s)^{-1}D_{s,x_s}J(C^{-1},x_s)uf(C^{-1}(x_s),u)\\
\\
=J_{-1}(C^{-1},x_s)^{-1}(I-P_{k,w})D_{s,x_s}J(C^{-1},x_s)\di\frac{J(C^{-1},x_s)wJ(C^{-1},x_s)}{\|J(C^{-1},x_s)\|^2}\\
\\
f(C^{-1}(x_s),\di\frac{J(C^{-1},x_s)wJ(C^{-1},x_s)}{\|J(C^{-1},x_s)\|^2})
\end{array}$$
Since $\di\frac{J(C^{-1},x_s)wJ(C^{-1},x_s)}{\|J(C^{-1},x_s)\|^2}=\di\frac{(x_s - e_{n+1})w(x_s - e_{n+1})}{||x_s - e_{n+1}||^2},$ the previous equation becomes
$$\begin{array}{ll}
Q_{k,u}uf(x,u)
=-J_{-1}(C^{-1},x_s)^{-1}Q_{k,w}^SwJ(C^{-1},x_s)f(C^{-1}(x_s),\di\frac{(x_s-e_{n+1})w(x_s-e_{n+1})}{\|(x_s-e_{n+1})\|^2}).\quad \blacksquare
\end{array}$$

Similarly, we have the following result for the remaining operators under the Cayley transformation.
\begin{theorem}
$$
-J_{-1}(C,x)Q_{k,u}^Sug(x_s,u)=Q_{k,w}wJ(C,x)g(C(x),\di\frac{(x+e_{n+1})w(x+e_{n+1})}{\|x+e_{n+1}\|^2}),
$$
where $u=\displaystyle\frac{(x+e_{n+1})w(x+e_{n+1})}{||x+e_{n+1}||^2},$
 $J(C,x)=\di\frac{x +e_{n+1}}{\|x+e_{n+1}\|^n}$ and $J_{-1}(C,x)=\di\frac{x+e_{n+1}}{\|x +e_{n+1}\|^{n+2}}$ is the conformal weight for the Cayley transformation.
\end{theorem}

As a consequence of two previous theorems we have the conformal invariance of equation $Q_k^S uf(x_s,u) = 0$:
\begin{theorem}
$Q_{k,u}uf(x,u)=0$ if and only if
$$
Q_{k,w}^SwJ(C^{-1},x_s)f(C^{-1}(x_s),\di\frac{(x_s-e_{n+1})w(x_s-e_{n+1})}{\|x_s-e_{n+1}\|^2})= 0
$$
and $Q_{k,u}^Sug(x_s,u) = 0 $ if and only if
$$
Q_{k,w}wJ(C,x)g(C(x),\di\frac{(x+e_{n+1})w(x+e_{n+1})}{\|x+e_{n+1}\|^2})= 0.
$$
\end{theorem}
\subsection{A kernel for the $Q_k^S$ operator }
Now consider the kernel in $ \mathbb{R}^n $
$$\begin{array}{ll}\di\frac{-1}{\omega_n c_k}w\di\frac{x-y}{\|x-y\|^n}Z_{k-1}(\di\frac{(x-y)w(x-y)}{\|x-y\|^2},v)v\\
\\
=\di\frac{-1}{\omega_n c_k}\di\frac{J(C^{-1},x_s)^{-1}uJ(C^{-1},x_s)^{-1}}{\|J(C^{-1},x_s)^{-1}\|^2}\\
\\
J(C^{-1},x_s)^{-1}\di\frac{x_s-y_s}{\|x_s-y_s\|^n}J(C^{-1},y_s)^{-1}Z_{k-1}(\di\frac{(x-y)w(x-y)}{\|x-y\|^2},v)v,
\end{array}$$ where $w=\di\frac{J(C^{-1},x_s)^{-1}uJ(C^{-1},x_s)^{-1}}{\|J(C^{-1},x_s)^{-1}\|^2}.$

Multiplying by $J(C^{-1},x_s)$ and applying the Cayley transformation to the above kernel, we obtain the kernel
\begin{eqnarray}\label{secondeq}
H_k^S(x-y,u,v):=
\di\frac{-1}{\omega_nc_k}u\di\frac{x_s-y_s}{\|x_s-y_s\|^n}J(C^{-1},y_s)^{-1}Z_{k-1}(au\tilde{a},v)v,
\end{eqnarray}
where $a=a(x_s,y_s)=\di\frac{J(C^{-1},x_s)^{-1}(x_s-y_s)J(C^{-1},y_s)^{-1}}{\|J(C^{-1},x_s)^{-1}\|\|x_s-y_s\|\|J(C^{-1},y_s)^{-1}\|}.$

This is a fundamental solution to $Q_k^Suf(x_s,u)= 0$ on $ \mathbb{S}^n,$ for $x_s,y_s \in \mathbb{S}^n.$

Similarly, we obtain that
\begin{eqnarray}\label{firsteq}
\di\frac{-1}{\omega_nc_k}uZ_{k-1}(u,\tilde{a}va)J(C^{-1},y_s)^{-1}\di\frac{x_s-y_s}{\|x_s-y_s\|^n}v
\end{eqnarray}
is a non trivial solution to $g(x_s, v)v Q_{k,r}^S= 0$.

We can see that the representations (\ref{secondeq}) and (\ref{firsteq}) are the same up to a reflection by
$$
\begin{array}{ll}
\di\frac{-1}{\omega_nc_k}uZ_{k-1}(u,\tilde{a}va)J(C^{-1},y_s)^{-1}\di\frac{x_s-y_s}{\|x_s-y_s\|^n}v\\
\\
=\di\frac{1}{\omega_nc_k}u\tilde{a}Z_k(au\tilde{a},v)aJ(C^{-1},y_s)^{-1}\di\frac{x_s-y_s}{\|x_s-y_s\|^n}v\\
\\
=\di\frac{-1}{\omega_nc_k}uJ(C^{-1},y_s)^{-1}\di\frac{x_s-y_s}{\|x_s-y_s\|^n}\di\frac{J(C^{-1},x_s)^{-1}}{\|J(C^{-1},x_s)^{-1}\|}
Z_k(au\tilde{a},v)\di\frac{J(C^{-1},x_s)^{-1}}{\|J(C^{-1},x_s)^{-1}\|}v\\
\\
=u\di\frac{J(C^{-1},x_s)^{-1}}{\|J(C^{-1},x_s)^{-1}\|}\di\frac{-1}{\omega_nc_k}\di\frac{x_s-y_s}{\|x_s-y_s\|^n}J(C^{-1},y_s)^{-1}
Z_k(au\tilde{a},v)\di\frac{J(C^{-1},x_s)^{-1}}{\|J(C^{-1},x_s)^{-1}\|}v.
\end{array}
$$

\subsection{Some basic integral formulas for the remaining operators on spheres}
In this section we will study some basic integral formulas related to the remaining operators on the sphere.

\begin{theorem} (Stokes' Theorem for the $n$-spherical Dirac operator $D_{s}$) \cite{LR}

Suppose $U_s$ is a domain on $\mathbb{S}^n$ and $f,g: U_s \times \mathbb{R}^{n} \to Cl_{n+1}$
are $C^1$, then for a subdomain $V_s$ of $U_s$, we have
$$\begin{array}{ll}
\di\int_{\partial V_s} g(x_s, u)n(x_s)f(x_s, u)d\Sigma(x_s)\\
\\
=\di\int_{V_s} (g(x_s, u)D_{s,x_s})f(x_s, u) +  g(x_s, u)(D_{s,x_s}f(x_s, u))dS(x_s) ,
\end{array}$$
where $dS(x_s)$ is the $n$-dimensional area measure on $V_s $, $d\Sigma(x_s)$ is the $n-1$-dimensional scalar
Lebesgue measure on $\partial V_s$ and $n(x_s)$ is the unit outward normal vector to $\partial V_s$ at $x_s$.
\end{theorem}

Applying the similar arguments to prove the Stokes' Theorem for $Q_k$ operators in section 4, we can obtain

\begin{theorem}(Stokes' Theorem for the $Q_k^S$ operator )

Let $U_s, V_s, \partial V_s$ be as in the previous Theorem. Then for $ f,g \in C^1(U_s \times \mathbb{R}^{n}, {\mathcal M}_{k})$,
we have
 version 1
$$\begin{array}{ll}
\di\int_{V_s}[(g(x_s,u)Q_{k,r}^S, f(x_s,u))_u+(g(x_s,u), Q_k^Sf(x_s,u))_u]dS(x_s)\\
\\
=\di\int_{\partial V_s}\left(g(x_s,u), (I-P_k)n(x_s)f(x_s,u)\right)_ud\Sigma(x)\\
\\
=\di\int_{\partial V_s}\left(g(x_s,u)n(x_s)(I-P_{k,r}), f(x_s,u)\right)_ud\Sigma(x).
\end{array}$$
Then for $f, g \in C^1(U_s \times \mathbb{R}^{n},$$\mathcal{M}_{k-1})$, we have version 2
$$\begin{array}{ll}
\di\int_{V_s}[(g(x_s,u)uQ_{k,r}^S, uf(x_s,u))_u+(g(x_s,u)u, Q_k^Suf(x_s,u))_u]dS(x_s)\\
\\
=\di\int_{\partial V_s}\left(g(x_s,u)u, (I-P_k)n(x_s)uf(x_s,u)\right)_ud\Sigma(x)\\
\\
=\di\int_{\partial V_s}\left(g(x_s,u)un(x_s)(I-P_{k,r}), uf(x_s,u)\right)_ud\Sigma(x).
\end{array}$$
\end{theorem}

\begin{remark}
Using the similar arguments to show the conformal invariance of Stokes' Theorem for the Rarita-Schwinger operators in \cite{LRV}, we obtain that Stokes' Theorem for the $Q_k$ operators is conformally invariant under the Cayley transformation and the inverse of the Cayley transformation.
\end{remark}
\begin{remark}
We also have the following fact
\begin{eqnarray}\label{7}
\di\int_{\partial V_s}\left(g(x_s,u)u, (I-P_k)n(x_s)uf(x_s,u)\right)_ud\Sigma(x)=\di\int_{\partial V_s}\left(g(x_s,u)u,n(x_s)uf(x_s,u)\right)_ud\Sigma(x)
\end{eqnarray}\end{remark}
\begin{theorem} (Borel-Pompeiu Theorem)
Suppose $U_s$, $V_s$ and $\partial V_s$ are stated as in Theorem $10$ and $y_s \in V_s.$
Then for $f \in C^1(U_s \times \mathbb{R}^{n},{\mathcal M}_{k-1}) $ we have
$$\begin{array}{ll}
u'f(y_s, u') =J(C^{-1},y_s)\di\int_{\partial V_s} (H_k^{S}(x_s-y_s, u, v), (I-P_{k})n(x_s) vf(x_s,v))_{v}d\Sigma(x_s) \\
\\
-J(C^{-1},y_s)\di\int_{V_s} (H_k^{S}(x_s-y_s, u, v),Q_k^Svf(x_s,v))_{v} dS(x_s)
\end{array}$$
where $u'=\di\frac{(y_s-e_{n+1})u(y_s-e_{n+1})}{\|y_s-e_{n+1}\|^2},$
$dS(x_s)$ is the $n$-dimensional area measure on $V_s \subset \mathbb{S}^n $, $n(x_s)$ and $d\Sigma(x_s)$ as in Theorem $10$.
\end{theorem}
{\bf Proof:} In the proof we use the representation
$$ \begin{array}{ll}H_k^S(x-y,u,v)=
\di\frac{-1}{\omega_nc_k}uZ_{k-1}(u,\tilde{a}va)J(C^{-1},y_s)^{-1}\di\frac{x_s-y_s}{\|x_s-y_s\|^n}v.
\end{array}$$

Let $B_s(y_s, \epsilon)$ be the ball centered at $y_s \in \mathbb{S}^n$ with radius $\epsilon$. We denote
$C^{-1}(B_s(y_s, \epsilon))$ by $B(y, r)$, and
$C^{-1}(\partial B_s(y_s, \epsilon))$ by $\partial B(y, r),$ where $y = C^{-1}(y_s)\in \mathbb{R}^n $ and $r$ is the radius of $B(y, r)$ in $\mathbb{R}^n$. Using the similar arguments in the proof of Theorem $2$, we only deal with
$$\begin{array}{ll}
\di\int_{\partial B_s(y_s, \epsilon)} (H_k^{S}(x_s-y_s, u, v),(I-P_k)n(x_s)vf(y_s,v))_{v} d\Sigma(x_s)\\
\\
=\di\int_{\partial B_s(y_s, \epsilon)}\int_{\mathbb{S}^{n-1}} H_k^{S}(x_s-y_s, u, v)(I-P_k)n(x_s)vf(y_s,v) ds(v)d\Sigma(x_s).\end{array}$$
Now applying (\ref{7}), the integral is equal to
$$\begin{array}{lll}
\di\int_{\partial B_s(y_s, \epsilon)}\int_{\mathbb{S}^{n-1}} H_k^{S}(x_s-y_s, u, v)n(x_s)vf(y_s,v) ds(v)d\Sigma(x_s)\\
\\
=\di\int_{\partial B_s(y_s, \epsilon)}\int_{\mathbb{S}^{n-1}} \di\frac{-1}{\omega_nc_k}uZ_{k-1}(u,\tilde{a}va)
J(C^{-1},y_s)^{-1}\di\frac{x_s-y_s}{\|x_s-y_s\|^n}v n(x_s)vf(y_s,v) ds(v)d\Sigma(x_s)\end{array}$$

Applying the inverse of the Cayley transformation to the previous integral, we have
$$\begin{array}{ll}
=\di\int_{\partial B(y, r)}\int_{\mathbb{S}^{n-1}} \di\frac{-1}{\omega_nc_k}u
Z_{k-1}(u,\di\frac{(x-y)w(x-y)}{\|x-y\|^2})
J(C^{-1},y_s)^{-1}J(C,y)^{-1}\di\frac{x-y}{\|x-y\|^n}J(C,x)^{-1}\\
\\
vJ(C,x)n(x)J(C,x)vf(C(y),\di\frac{J(C,y)wJ(C,y)}{\|J(C,y)\|^2}) ds(v)d\sigma(x),
\end{array}$$
where $v=\di\frac{J(C,y)wJ(C,y)}{\|J(C,y)\|^2}.$

Now if we replace $v$ with $\di\frac{J(C,y)wJ(C,y)}{\|J(C,y)\|^2}$ in the previous integral and we also set $J(C,x)=(J(C,x)-J(C,y))+J(C,y)$, but $J(C,x)-J(C,y)$ tends to zero as $x$ approaches $y$. Thus the previous integral
can be replaced by
$$\begin{array}{ll}
=\di\int_{\partial B(y, r)}\int_{\mathbb{S}^{n-1}} \di\frac{-1}{\omega_nc_k}u
Z_{k-1}(u,\di\frac{(x-y)w(x-y)}{\|x-y\|^2})
\di\frac{x-y}{\|x-y\|^n}J(C,y)^{-1}\\
\\
\di\frac{J(C,y)wJ(C,y)}{\|J(C,y)\|^2}J(C,y)n(x)J(C,y)\di\frac{J(C,y)wJ(C,y)}{\|J(C,y)\|^2}f(C(y),\di\frac{J(C,y)wJ(C,y)}{\|J(C,y)\|^2}) ds(w)d\sigma(x)\\
\\
=\di\int_{\partial B(y, r)}\int_{\mathbb{S}^{n-1}} \di\frac{-1}{\omega_nc_k}u
Z_{k-1}(u,\di\frac{(x-y)w(x-y)}{\|x-y\|^2})
\di\frac{x-y}{\|x-y\|^n}
wn(x)\\
\\
wJ(C,y)f(C(y),\di\frac{J(C,y)wJ(C,y)}{\|J(C,y)\|^2}) ds(w)d\sigma(x)\\
\\
=\di\int_{\partial B(y, r)}\int_{\mathbb{S}^{n-1}} \di\frac{1}{\omega_nc_k}\di\frac{1}{r^{n-1}}u
Z_{k-1}(u,\di\frac{(x-y)w(x-y)}{\|x-y\|^2})
\di\frac{x-y}{\|x-y\|}
w\di\frac{x-y}{\|x-y\|}\\
\\
wJ(C,y)f(C(y),\di\frac{J(C,y)wJ(C,y)}{\|J(C,y)\|^2}) ds(w)d\sigma(x).
\end{array}$$

Using Lemma $5$ in \cite{DLRV}, the previous integral becomes
\begin{eqnarray}\label{8}
&&\di\int_{\mathbb{S}^{n-1}}uZ_{k-1}(u,w)wwJ(C,y)f(C(y),\di\frac{J(C,y)wJ(C,y)}{\|J(C,y)\|^2})ds(w)\nonumber\\
\nonumber\\
&&=-\di\int_{\mathbb{S}^{n-1}}uZ_{k-1}(u,w)J(C,y)f(C(y),\di\frac{J(C,y)wJ(C,y)}{\|J(C,y)\|^2})ds(w)\nonumber\\
\nonumber\\
&&=-uJ(C,y)f(C(y),\di\frac{J(C,y)uJ(C,y)}{\|J(C,y)\|^2}).
 \end{eqnarray}
If we set $u'=\di\frac{J(C,y)uJ(C,y)}{\|J(C,y)\|^2}=\di\frac{J(C^{-1},y_s)^{-1}uJ(C^{-1},y_s)^{-1}}{\|J(C^{-1},y_s)^{-1}\|^2}=\di\frac{(y_s-e_{n+1})u(y_s-e_{n+1})}{\|y_s-e_{n+1}\|^2},$ then\\ $uJ(C,y)=uJ(C^{-1},y_s)^{-1}=J(C^{-1},y_s)\|J(C^{-1},y_s)^{-1}\|^2u'.$ Now if we multiply the both sides of equation (\ref{8}) by
$\di\frac{J(C^{-1},y_s)^{-1}}{\|J(C^{-1},y_s)^{-1}\|^2}=-J(C^{-1},y_s)$, then we obtain
$$\begin{array}{ll}
J(C^{-1},y_s)\di\int_{\mathbb{S}^{n-1}}uZ_{k-1}(u,w)J(C,y)f(C(y),\di\frac{J(C,y)wJ(C,y)}{\|J(C,y)\|^2})ds(w)\\
\\
=u'f(C(y),u')=u'f(y_s,u').\quad \blacksquare\end{array}$$

\begin{corollary}

Let $\psi$ be a function in $C^\infty (V_s,\mathcal{M}_{k-1})$ and supp $f\subset V_s$. Then
$$
u'\psi(y_s, u') = - J(C^{-1},y_s)\int_{V_s} (H_k^{S}(x_s-y_s, u, v),Q_k^Sv\psi(x_s,v))_{v} dS(x_s),
$$
where $u'=\di\frac{(y_s-e_{n+1})u(y_s-e_{n+1})}{\|y_s-e_{n+1}\|^2}.$
\end{corollary}

\begin{corollary} (Cauchy Integral Formula for $Q_k^S$ operators)

If $Q_k^Svf(x_s, v) = 0$, then for $y_s \in V_s$ we have
\begin{eqnarray*}
u'f(y_s, u') &=&J(C^{-1},y_s)\int_{\partial V_s} (H_k^{S}(x_s-y_s, u, v), (I-P_k)n(x_s)v f(x_s,v))_{v} d\Sigma(x_s) \\
&=&J(C^{-1},y_s)\int_{\partial V_s} (H_k^{S}(x_s-y_s, u, v) n(x_s)(I-P_{k,r}),v f(x_s,v))_{v} d\Sigma(x_s),
\end{eqnarray*}
where $u'=\di\frac{(y_s-e_{n+1})u(y_s-e_{n+1})}{\|y_s-e_{n+1}\|^2}.$
\end{corollary}

\begin{remark}
By factoring out $\mathbb{S}^n$ by the group $\mathbb{Z}_2=\{\pm 1\}$ we obtain real projective space, $\mathbb{R}P^n$. Using the similar arguments to obtain the results for Rarita-Schwinger operators on real projective space in \cite{LRV}, we can easily extends the similar results for $Q_k$ operators to real projective space.
\end{remark}

%%%%%%%%%
Junxia Li  \qquad  Email: jxl004@uark.edu\\
John Ryan   \qquad Email:  jryan@uark.edu
\end{document}